\documentclass[10pt,twoside]{article}
\usepackage{float,latexsym,graphicx,shortvrb,fancyhdr}
\usepackage{amsmath,amsfonts,bm,mathtools,caption, marvosym}
\usepackage[a4paper]{geometry}
\fancypagestyle{firststyle}
{
   
   \fancyhf{}
   \fancyhead[L]{}
   \fancyhead[R]{}
   \fancyfoot[L]{ \footnotesize Y. Chen $\cdot$ D.Y. Gao (\Letter)  \\ School of Science, Information Technology and Engineering, Federation University Australia, Victoria, Australia,\\ e-mail: yi.chen@federation.edu.au\\
   \vspace{1em}
   D.Y. Gao\\
   e-mail: d.gao@federation.edu.au
   }
   }

\makeatletter
\renewcommand{\maketitle}{\bgroup\setlength{\parindent}{0pt}
\begin{flushleft}
  \textbf{\@title}

  \@author
\end{flushleft}\egroup
}
\makeatother

\usepackage{amsmath,amsfonts,bm,mathtools,graphicx,latexsym,shortvrb}
\usepackage{xcolor}
\usepackage{hyperref}
\usepackage{multirow}
\usepackage{colortbl}

\newtheorem{theorem}{Theorem}
\newtheorem{lemma}[theorem]{Lemma}
\newtheorem{proposition}[theorem]{Proposition}
\newcommand{\proof}{\par \noindent \bf Proof: \hspace{0mm} \rm }
\newcommand{\qed}{\hspace*{\fill} $\Box$ \vspace{2ex}}

\def\btheo{\begin{theorem}}
\def\etheo{\end{theorem}}
\def\bprop{\begin{proposition}}
\def\eprop{\end{proposition}}
\def\bexam{\begin{example}}
\def\eexam{\end{example}}
\def\bdefi{\begin{definition}}
\def\edefi{\end{definition}}
\def\blemm{\begin{lemma}}
\def\elemm{\end{lemma}}

\newcommand{\mbb}{\mathbb}
\newcommand{\mcal}{\mathcal}

\newcommand{\sta}{\textrm{sta}}
\def\inv{{-1}}

\def\det{\textrm{det}}
\def\diag{\textrm{diag}}

\def\exp{\textrm{exp}}

\def\rank{\textrm{rank}}

\def\trace{\textrm{trace}}
\def\non{\nonumber}
\def\wand{\textrm{ and }}
\def\wotherwise{\textrm{otherwise}}

\def\[#1\]{\begin{align}#1\end{align}}
\def\bcase{\begin{cases}}
\def\ecase{\end{cases}}
\def\bpmat{\begin{pmatrix}}
\def\epmat{\end{pmatrix}}
\def\bbmat{\begin{bmatrix}}
\def\ebmat{\end{bmatrix}}
\def\beqn{\begin{eqnarray}}
\def\eeqn{\end{eqnarray}}
\def\beqnx{\begin{eqnarray*}}
\def\eeqnx{\end{eqnarray*}}
\def\beq{\begin{equation}}
\def\eeq{\end{equation}}
\def\bitem{\begin{itemize}}
\def\eitem{\end{itemize}}
\def\btheo{\begin{theorem}}
\def\etheo{\end{theorem}}
\def\bblock{\begin{block}}
\def\eblock{\end{block}}
\def\benum{\begin{enumerate}}
\def\eenum{\end{enumerate}}

\def\barx{\bar{x}}
\def\bxx{\bar{\bm x}}

\def\hatf{\hat{f}}
\def\hff{\hat{\bm f}}

\def\xx{\bm x}

\def\ff{\bm f}
\def\ee{\bm e}

\def\cc{\bm c}
\def\dd{\bm d}

\def\II{\bm I}
\def\AA{\bm A}
\def\BB{\bm B}

\def\DD{\bm D}
\def\EE{\bm E}
\def\FF{\bm F}
\def\GG{\bm G}
\def\KK{\bm K}
\def\LL{\bm L}
\def\PP{\bm P}
\def\QQ{\bm Q}
\def\RR{\bm R}

\def\UU{\bm U}
\def\XX{\bm X}

\def\a{\alpha}

\def\b{\beta}

\def\gcc{\bm\chi}

\def\d{\delta}

\def\gzz{\bm\zeta}

\def\gthth{\bm\theta}

\def\gvthth{\bm\vth}

\def\l{\lambda}
\def\gll{\bm\lambda}

\def\s{\sigma}
\def\gss{\bm\sigma}
\def\t{\tau}
\def\gtt{\bm\tau}

\def\x{\xi}
\def\gxx{\bm\xi}

\def\gyy{\bm\eta}

\def\vP{\varPi}

\def\vX{\varXi}

\def\bgzz{\bar{\bm\zeta}}
\def\bgs{\bar{\s}}

\def\bgt{\bar{\t}}

\def\bbR{\mathbb{R}}

\def\calE{\mathcal{E}}

\def\calP{\mathcal{P}}

\def\bbS{\mathbb{S}}
\def\calS{\mathcal{S}}

\def\calX{\mathcal{X}}

\def\lpa{\left(}
\def\rpa{\right)}
\def\lbk{\left[}
\def\rbk{\right]}
\def\lbc{\left\{}
\def\rbc{\right\}}

\def\half{\frac{1}{2}}

\renewcommand{\x}{\chi}
\renewcommand{\gxx}{\mbox{\boldmath$\chi$}}
\renewcommand{\gcc}{\mbox{\boldmath$\xi$}}
\renewcommand{\gvthth}{\mbox{\boldmath$\omega$}}
\begin{document}

 \thispagestyle{firststyle}

\title{\bf\Large Global Solutions to Nonconvex Optimization of 4th-Order Polynomial and Log-Sum-Exp Functions}
\maketitle

\vspace{3em}

\noindent{{\bf Yi Chen $\cdot$ David Y Gao}}

\vspace{3em}

\noindent{\textbf{Abstract}}~
This paper presents a canonical dual approach for solving a nonconvex global optimization problem
governed by a sum of fourth-order polynomial and a log-sum-exp function.
Such a problem arises extensively in  engineering and sciences. Based on the canonical duality-triality  theory,
 this nonconvex problem is transformed to an equivalent dual problem, which can be solved  easily
 under certain conditions.
We proved that  both  global minimizer and the biggest local extrema  of the primal problem can be obtained analytically from the canonical dual solutions.
  As two special cases,  a quartic  polynomial minimization   and  a minimax problem are discussed.
   Existence conditions   are derived, which can be used to classify easy and relative hard instances.
  Applications are illustrated by several  nonconvex and nonsmooth  examples.

\vspace{2em}


\vspace{2em}


\vspace{3em}
\textheight=24cm

\section{Introduction}\label{se:intro}
In this paper, we are interested in the following nonconvex global optimization problem:
\begin{equation}\label{p:primal}
(\mathcal{P}):~~~~\min_{\xx\in\bbR^n} \vP(\bm x)= \frac{1}{2}\bm x^T \AA\bm x-\bm f^T\bm x +T(\bm x)+W(\bm x)
\end{equation}
in which $\AA\in\mathbb{R}^{n\times n}$ is a symmetric matrix, $\bm f\in\mathbb{R}^n$, and the log-sum-exp function $T(\bm x)$ and the fourth-order polynomial function $W(\bm x)$ are defined as
\begin{eqnarray}
&&T(\bm x)=\frac{1}{\beta}\log\lbk1+\sum_{i=1}^p\exp\lpa\beta\lpa\frac{1}{2}\bm x^T \QQ_i
\xx+d_i\rpa\rpa\rbk,\label{eq:funLse}\\
&&W(\bm x)= \sum_{i=1}^r\frac{\alpha_i}{2}\lpa\frac{1}{2}\bm x^T \BB_i\bm x+c_i\rpa^2,\label{eq:funW}
\end{eqnarray}
where $ \QQ_i\in\mathbb{R}^{n\times n}$ and $ \BB_i\in\mathbb{R}^{n\times n}$ are symmetric matrices, $d_i$ and $c_i$ are any real numbers, and $\alpha_i$ and $\beta$ are positive real numbers.

The quartic function $W(\bm x)$ is the so-called double-well potential if $ \BB_i \succeq 0, \; c_i < 0$ for $i = 1, \dots, r$. This function has extensive applications in
 mathematical physics, for example, in \cite{gao00buckling} it was used to model post-buckling of beams. While  $ T(\bm x)$
 is one of the fundamental functions in engineering sciences,  which arises broadly in regions including plasticity theory \cite{strang79minimax},
  nonsmooth variational problems \cite{gao98minimax}, structural optimization problems \cite{banichuk76minimax},   information theory \cite{chiang04gp},  network communication systems \cite{chiang05GP,boyd07GP,chiang07power},
and robot manipulator designing \cite{abdi10jacobian,abdi11jacobian,roberts07robot}.
In numerical analysis,  the function  $ T(\bm x)$  is often used to deal with minimax problems
 \cite{pee11LSE,polak97optimization,polak03minimax,royset04minimax}.
 Due to the nonconvexity,  it has been proven that directly solving $(\mathcal{P})$ is very difficult.
 The main difficulty is the global optimality condition.
 Although the traditional Lagrangian and modern Fenchel-Moreau-Rockaffelar  duality methods can be used to produce lower bound approaches,
  these methods suffer for having duality gap problems.
 Therefore, a perfect duality theory (with zero duality gap)  should be fundamentally important not only for solving nonconvex problems, but also for understanding
 nonconvex phenomena in complex systems.

Canonical duality theory was originally developed from infinite dimensional  nonconvex systems   \cite{Gao00nonsmooth}.
This potentially powerful  theory has been applied successfully for solving many challenging  problems  in global optimization and nonconvex/nonsmooth analysis,
such as quadratic problems \cite{Fang08qp01,Gao04quadratic,Gao10boxIP,Gao09Lagrangian}, polynomial optimization \cite{Gao06polynomial},
transportation problems \cite{gao10fixedcost}, location problems \cite{Gao10sensor}, and integer programming  problems \cite{wang11maxcut}.
 Recently, after an open problem left on the triality theory been solved \cite{Gao-Wu-triality12}
 and some efficient algorithms been developed  \cite{gao11qp}, the canonical duality-triality theory is considered as a breaking-through theory in the multi-disciplinary
 fields of mathematical physics and global optimization.

The purpose of this paper is to apply the canonical duality theory to solve globally the nonconvex optimisation problem  $(\mathcal{P})$.
The rest of this paper is arranged as follows.
 We first  show in Section \ref{se:cano} how a canonical dual problem can be  constructed by the standard canonical dual transformation.
 Then in Section \ref{se:triality}, the  triality theory is presented and proved, which shows
 that the nonconvex problem $(\mathcal{P})$ is actually equivalent to a convex problem if its canonical dual has a critical point in the positive semidefinite region.
  Correspondingly, the global solution of the primal problem can be obtained analytically from this critical point.
Two specific problems are discussed in Section \ref{se:specific}, one is a nonconvex quartic  polynomial minimization problem and one is a nonsmooth
 minimax problem.   Examples are provided in Section \ref{se:exs} to illustrate the canonical duality theory.
Some conclusions are given in Section \ref{se:concl}.

\section{Canonical dual problem and analytical solutions}\label{se:cano}
 Following  standard procedure of the canonical dual transformation,
we need to introduce  a so-called geometrically admissible operator
$$
\gcc(\xx)=\bpmat \gxx(\xx)\\\gyy(\xx) \epmat=\bpmat
\lbc\frac{1}{2}\bm x^T \QQ_i\bm x\rbc_{i=1}^p\vspace{.2cm}\\
\lbc\frac{1}{2}\bm x^T \BB_i\bm x\rbc_{i=1}^r
\epmat
~:~\mathbb{R}^n\rightarrow\mcal{E}_a\subseteq\mathbb{R}^m
$$
where $m=p+r$.
 Generally speaking,  for any given $m$ and symmetrical matrices  $\QQ_i$ and $\BB_i$,
  the range  $\calE_a$  may not be convex.   However,  if $m=2$,
  the range  $\calE_a$ is always convex.
  Suppose $p=r=1$. As proved in \cite{Boyd04convex}, we will have $\calE_a=\{(\half \trace(\QQ_1\XX),\half\trace(\BB_1\XX))~|~\XX\in\bbS_+^n\}$, where the right side of the equation is the range
   of all positive semidefinite matrix under a linear transformation and thus it is a convex set. In this paper,
    we assume that $\calE_a$ is a convex set.
Therefore, a canonical function can be  defined on  $\calE_a$:
$$
V(\gcc)=V_1(\gxx)+V_2(\gyy)
$$
where
\begin{eqnarray}
&&V_1(\gxx)= \frac{1}{\beta}\log\lbk1+\sum_{i=1}^p\exp\lpa\beta (\x_i+d_i)\rpa\rbk,\non\\
&&V_2(\bm\eta)=\sum_{i=1}^r\frac{\alpha_i}{2}(\eta_i+c_i)^2.\non
\end{eqnarray}
Here $\x_i$ denotes the $i$th component of $\gxx$, and $\eta_i$ denotes the $i$th component of $\bm\eta$. Since $V_1(\gxx)$ and $V_2(\bm\eta)$ are convex,
$V(\gcc)$ is a convex function. Thus, by Legendre transformation, we have the following relationship
\beq\label{eq:legendre}
V(\gcc)+V^*(\gzz)=\gcc^T\gzz,
\eeq
where
$$
\gzz=\bpmat \gtt\\\gss \epmat
=\bpmat  \nabla V_1(\gxx)\\\nabla V_2(\bm\eta) \epmat
=\bpmat \lbc\frac{\exp(\beta(\x_i+d_i))}{1+\sum_{k=1}^p\exp(\beta(\x_k+d_k))}\rbc_{i=1}^p\\
 \lbc\alpha_i(\eta_i+c_i)\rbc_{i=1}^r\epmat ~:~\mcal{E}_a\rightarrow \mcal{E}_a^*\subseteq\bbR^m
$$
and $V^*(\gzz)$ is the conjugate function of $V(\gcc)$, defined as
\[
V^*(\gzz)=V_1^*(\bm\tau)+V_2^*(\bm\sigma)
\]
with
\begin{eqnarray}
&&V_1^*(\bm\tau)=\frac{1}{\beta}\lbk\sum_{i=1}^p \tau_i\log(\tau_i)+(1-\sum_{i=1}^p\t_i)\log(1-\sum_{i=1}^p\t_i)\rbk-\dd^T\gtt,\label{eq:funV1star}\non\\
&&V_2^*(\bm\sigma)=\sum_{i=1}^r\frac{1}{2\alpha_i}\sigma_i^2-\cc^T\gss,\label{eq:funV2star}\non
\end{eqnarray}
where $\dd=\{d_i\}$ and $\cc=\{c_i\}$.
It can be verified that the Jacobian of the map $\gzz$ from $\calE_a$ to $\calE_a^*$ is not zero, which indicates that the map is invertible.
 Since $\calE_a$ is  assumed to be convex, the image $\calE_a^*$ is also convex. Moreover, from the definition of $\gtt$,
it is required that $\gtt \in   \{\gtt~:~ \gtt> 0, \;\; \gtt^T\ee<1\}$. In this paper, we denote $\ee$ as a vector whose components are all ones. 

The so-called generalized total complementary function
$\vX : \mbb{R}^n\times \mcal{E}_a^*\rightarrow\mbb{R}$ can then be defined as
\[\label{eq:totalcompf}
\vX(\bm x,\bm\zeta) 
&=\frac{1}{2}\bm x^T\AA\bm x-\bm f^T\bm x+\gcc^T\gzz-V^*(\gzz)\nonumber\\
&=\frac{1}{2}\bm x^T \GG_a(\gzz)\bm x-\ff^T\bm x-V_1^*(\bm\tau)-V_2^*(\bm\sigma),
\]
where
$$
\GG_a(\gzz)= \AA+\sum_{i=1}^p\tau_i \QQ_i+\sum_{i=1}^r\sigma_i \BB_i.
$$

From the generalized total complementary function,   the canonical dual function $\vP^d(\bm\zeta)$  can be obtained by
\begin{equation}\label{eq:dualfun0}
\vP^d(\bm\zeta)=\sta\lbc\vX(\bm x,\bm\zeta)~|~\bm x\in\mbb{R}^n\rbc, ~~ \gzz\in\calE_a^*
\end{equation}
where the notation $\sta\{\cdot\}$ represents the task of finding stationary points of $\vX(\bm x,\bm\zeta)$ with respect to $\bm x$. Notice that for any given $\bm\zeta$, the total complementary function $\vX(\bm x,\bm\zeta)$ is a quadratic function of $\bm x$ and its stationary points are the solutions of the following equation system
\begin{equation}\label{eq:partialXi}
\nabla_{\bm x}\vX(\bm x,\bm\zeta)= \GG_a(\gzz)\bm x-\ff=0.
\end{equation}
If $\det(\GG_a(\gzz))\neq 0$ for a given $\gzz$, $\bm x$ can be solved analytically and uniquely as
$\bm x= \GG_a(\gzz)^\inv\ff$ and the canonical dual function $\vP^d(\bm\zeta)$ can then be written explicitly as
\[\label{eq:dualfun}
\vP^d(\bm\zeta)=&-\frac{1}{2}\ff^T \GG_a(\gzz)^\inv\ff-V_1^*(\bm\tau)-V_2^*(\bm\sigma).
\]

Let
$$
\mcal{S}_a=\lbc\bm\zeta~|~\bm\zeta\in\mcal{E}_a^*,~ \det(\GG_a(\gzz))\neq0 \rbc.
$$
The canonical dual problem can be proposed as the following. 
\begin{equation}\label{p:dual}
(\mcal{P}^d):~~~~\sta\lbc\vP^d(\bm\zeta)~|~\bm\zeta\in\mcal{S}_a\rbc.
\end{equation}

\begin{theorem}\label{th:AnalSolu}
{\rm (\textbf{Analytical Solution and Complementary-Dual Principle} \cite{Gao00duality,Gao-Wu-triality12})}
The problem ($\calP^d$) is canonically dual to the problem ($\calP$) in the sense that if $\bar{\bm\zeta}\in\mcal{S}_a$ is a critical point of $\vP^d(\bm\zeta)$, then
\beq\label{eq:solvedx}
\bar{\bm x}=\GG_a(\bgzz)^\inv\ff
\eeq
is a critical point of $\vP(\bm x)$, the pair $(\bxx,\bgzz)$ is a critical point of $\vX(\xx,\gzz)$, and we have
\beq\label{eq:nogap}
\vP(\bar{\bm x})=\vX(\bxx,\bgzz)=\vP^d(\bar{\bm\zeta}).
\eeq
\end{theorem}
The proof of this theorem  is  analogous with that in \cite{Gao00duality}.
Theorem \ref{th:AnalSolu} shows that  there is no duality gap between the primal problem ($\calP$) and the canonical dual problem ($\calP^d$).

\section{Triality theory}\label{se:triality}
In this section we study the global optimality conditions for the critical solutions of   the primal and dual problems.
We need to define  the following two subsets  of the dual space $\mcal{S}_a$:
\begin{eqnarray*}
&&\mcal{S}_a^+= \lbc\bm\zeta\in\mcal{S}_a~|~\GG_a\succ 0\rbc,\label{eq:Saplus}\\
&&\mcal{S}_a^-= \lbc\bm\zeta\in\mcal{S}_a~|~\GG_a\prec 0\rbc.\label{eq:Saminus}
\end{eqnarray*}
It is easy to  prove that both $\calS_a^+$ and $\calS_a^-$ are convex sets.
 Here we  use $\GG_a$ for short to denote $\GG_a(\gzz)$.

For convenience, we here give the first and second derivatives of functions $\vP(\xx)$ and $\vP^d(\gzz)$:
\begin{eqnarray}
&& \nabla \vP(\bm x)= \GG_a\bm x-\ff,\label{eq:der1Pi}\\
&& \nabla^2\vP(\bm x)= \GG_a+\FF \DD \FF^T,\label{eq:der2Pi}\\
&& \nabla \vP^d(\bm\zeta)=
\left(\begin{array}{l}
\lbc\frac{1}{2}\ff^T \GG_a^\inv \QQ_i \GG_a^\inv\ff+d_i-\frac{1}{\beta}\log\frac{\t_i}{1-\sum_{i=1}^p\t_i}\rbc_{i=1}^p\\
\lbc\frac{1}{2}\ff^T \GG_a^\inv \BB_i \GG_a^\inv\ff+c_i-\frac{\sigma_i}{\alpha_i}\rbc_{i=1}^r
\end{array}\right),\label{eq:der1Pid}\\
&& \nabla^2\vP^d(\bm\zeta)=- \FF^T \GG_a^\inv \FF-\DD^\inv,\label{eq:der2Pid}
\end{eqnarray}
where $\FF\in\mbb{R}^{n\times m}$ and $\DD\in\mbb{R}^{m\times m}$ are defined as
\begin{eqnarray*}
&& \FF=
\begin{bmatrix}
 \QQ_1\xx,\ldots, \QQ_p\xx, \BB_1\xx,\ldots, \BB_r\xx
\end{bmatrix},\label{eq:matrixF}\\
&& \bm \DD=
\begin{bmatrix}
\beta\lpa\diag(\bm\tau)-\bm\tau\bm\tau^T\rpa&0\\
0 &\diag(\bm\alpha)
\end{bmatrix} , \label{eq:matrixD}
\end{eqnarray*}
in which,  $\xx = \GG_a^\inv\ff$. By the fact that 
  $(\diag(\gtt)-\gtt\gtt^T)^\inv=\diag(\gtt)^\inv+\ee\ee^T/(1-\gtt^T\ee) \succ \;\; \forall \gtt \in \{\gtt~:~ \gtt > 0, \;\; \gtt^T\ee<1\}$, 
the matrix $\DD^\inv$ is positive definite. 

The following two lemmas are needed. Their proofs are omitted, which are analogous with that in \cite{Gao-Wu-triality12}.
\begin{lemma}\label{lm:hessianPi}
Suppose that $ m<n$, $\bar{\bm\zeta}\in\mcal{S}_a^-$ is a critical point and a local minimizer of $\vP^d(\bm\zeta)$, and $\bar{\bm x}=\GG_a^\inv\ff$. Then, there exists a matrix $ \LL\in\mbb{R}^{n\times m}$ with $\rank( \LL)=m$ such that
\begin{equation}\label{eq:LhessianPiL}
 \LL^T\nabla^2\vP(\bar{\bm x}) \LL\succeq0.
\end{equation}
\end{lemma}

\begin{lemma}\label{lm:hessPid}
Suppose that $ m>n$, $\bar{\bm\zeta}\in\mcal{S}_a^-$ is a critical point of $\vP^d(\bm\zeta)$, and $\bar{\bm x}=\GG_a^\inv\ff$ is a local minimizer of $\vP(\bm x)$. Then, there exists a matrix $ \PP\in\mbb{R}^{m\times n}$ with $\rank( \PP)=n$ such that
\begin{equation}\label{eq:lm2PhessPidP}
 \PP^T\nabla^2\vP^d(\bar{\bm\zeta}) \PP\succeq0.
\end{equation}
\end{lemma}
Accordingly, two sets are defined, which are generated by shifting the column spaces of $\LL$ and $\PP$:
\begin{eqnarray*}
&&\mcal{X}_{ \LL}= \lbc\bm x\in\mbb{R}^n~|~ \bar{\bm x}+\LL \gthth,~ \gthth\in\mbb{R}^m\rbc,\label{eq:setXL}\\
&&\mcal{S}_{ \PP}= \lbc\bm\zeta\in\mbb{R}^m~|~ \bar{\bm\zeta}+\PP\gvthth,~\gvthth\in\mbb{R}^n\rbc.\label{eq:setSP}
\end{eqnarray*}

Now, we give the main result of this paper, triality theorem, which illustrates the relationships between the primal and canonical dual problems on global and local solutions.
\begin{theorem}\label{th:triality}
{\rm(\textbf{Triality Theorem})}  Suppose that $\bar{\bm\zeta}$ is a critical point of $\vP^d(\bm\zeta)$, and $\bar{\bm x}=\GG_a(\bgzz)^\inv\ff$.
\benum
\item If $\bar{\bm\zeta}\in\mcal{S}_a^+$, then the canonical min-max duality holds in the form of
\begin{equation}\label{eq:th2minmax}
\vP(\bar{\bm x})   =\min_{\bm x\in\mbb{R}^n} \vP(\bm x)=\max_{\bm\zeta\in\mcal{S}_a^+} \vP^d(\bm\zeta)=\vP^d(\bar{\bm\zeta}).
\end{equation}
\item If $\bar{\bm\zeta}\in\mcal{S}_a^-$, the double-max duality holds in the form that if $\bxx$ is a local maximizer of $\vP(\bm x)$ or $\bgzz$ is a local maximizer of $\vP^d(\bm\zeta)$,
 we have
    \begin{equation}\label{eq:th2maxmax}
    \vP(\bar{\bm x})   =\max_{\bm x\in\mcal{X}_0} \vP(\bm x)=\max_{\bm\zeta\in\mcal{S}_0} \vP^d(\bm\zeta)=\vP^d(\bar{\bm\zeta})
    \end{equation}
for some neighborhood\footnote{We use the same definition of the neighborhood  as defined in \cite{gao-amma03} (Note 1 on page 306), i.e.,
a subset ${\cal X}_0$ is  said to be the neighborhood of the critical point
$\bar{\xx}$ if $\bar{\xx}$ is the only critical point in ${\cal X}_0$.}
 $\mcal{X}_0\times\mcal{S}_0\subset\mbb{R}^n\times\mcal{S}_a^-$ of $(\bar{\bm x},\bar{\bm\zeta})$.
\item If $\bgzz\in\calS_a^-$, then the double-min duality holds conditionally as:
\benum
\item Given $ m=n$, if $\bxx$ is a local minimizer of $\vP(\bm x)$ or $\bgzz$ is a local minimizer of $\vP^d(\bm\zeta)$, we have
    \begin{equation}\label{eq:th2minmin}
    \vP(\bar{\bm x})   =\min_{\bm x\in\mcal{X}_0} \vP(\bm x)=\min_{\bm\zeta\in\mcal{S}_0} \vP^d(\bm\zeta)=\vP^d(\bar{\bm\zeta})
    \end{equation}
for some neighborhood $\mcal{X}_0\times\mcal{S}_0\subset\mbb{R}^n\times\mcal{S}_a^-$ of  $(\bar{\bm x},\bar{\bm\zeta})$.
\item  Given $ m<n$, if $\bar{\bm\zeta}$ is a local minimizer of $\vP^d(\bm\zeta)$, $\bar{\bm x}$ will be a saddle point of $\vP(\bm x)$ and there exists a neighborhood $(\mcal{X}_0\cap\calX_{\LL})\times\mcal{S}_0\subset\mbb{R}^n\times\mcal{S}_a^-$ of $(\bar{\bm x},\bar{\bm\zeta})$ such that
    \begin{equation}\label{eq:th2minmln}
\vP(\bar{\bm x})   =\min_{\bm x\in\mcal{X}_0\cap\mcal{X}_{ \LL}} \vP(\bm x)=\min_{\bm\zeta\in\mcal{S}_0} \vP^d(\bm\zeta)=\vP^d(\bar{\bm\zeta});
\end{equation}
\item Given $ m>n$, if $\bar{\bm x}$ is a local minimizer of $\vP(\bm x)$, $\bar{\bm\zeta}$ will be a saddle point of $\vP^d(\bm\zeta)$ and there exists a neighborhood $\mcal{X}_0\times(\mcal{S}_0)\cap\calS_P\subset\mbb{R}^n\times\mcal{S}_a^-$ of $(\bar{\bm x},\bar{\bm\zeta})$ such that
\begin{equation}\label{eq:th2minmgn}
\vP(\bar{\bm x})=\min_{\bm x\in\mcal{X}_0} \vP(\bm x)=\min_{\bm\zeta\in\mcal{S}_0\cap\mcal{S}_{ P}} \vP^d(\bm\zeta)=\vP^d(\bar{\bm\zeta}).
\end{equation}
\eenum
\eenum
\end{theorem}
\proof
\benum
\item  Since $ \GG_a(\bgzz)\succ0$ when $\bar{\bm\zeta}\in\mcal{S}_a^+$ and $ \DD\succ0$, the Hessian of the dual function is negative definitive, i.e. $\nabla^2\vP^d(\bm\zeta)\prec0$,  which implies that $\vP^d(\bm\zeta)$ is strictly concave over $\mcal{S}_a^+$. Thus, we have
\[\label{eq:th2Pidmax}
\vP^d(\bgzz)=\max_{\gzz\in\calS_a^+} \vP^d(\gzz).
\]
Whereas, the generalized total complementary function $\vX(\bm x,\bar{\bm\zeta})$ is a convex function with respect to $\xx$ in $\mbb{R}^n$, which, plus the fact that $\bar{\bm x}$ is a critical point of $\vX(\bm x,\bar{\bm\zeta})$, implies that we have $\vX(\bm x,\bar{\bm\zeta})\geq\vX(\bar{\bm x},\bar{\bm\zeta})$ for any $\bm x\in\mbb{R}^n$. From the Fenchel's inequality, it is true that $
    \vX(\bm x,\bm\zeta)\leq\vP(\bm x), ~~\forall(\bm x,\bm\zeta)\in\mbb{R}^n\times\mcal{S}_a$.
    Therefore, for any $\bm x\in\mbb{R}^n$, we have
    \begin{equation}\label{eq:th2PiXiPi}
    \vP(\bm x)\geq\vX(\bm x,\bar{\bm\zeta})\geq\vX(\bar{\bm x},\bar{\bm\zeta})=\vP(\bar{\bm x}).
    \end{equation}
 Thus, the equation (\ref{eq:th2minmax}) is true.

\item Suppose $\bar{\bm\zeta}$ is a local maximizer of $\vP^d(\bm\zeta)$ in $\mcal{S}_a^-$. Then we have $\nabla^2\vP^d(\bar{\bm\zeta})=- \FF^T \GG_a^\inv \FF- \DD^{-1}\preceq0$ and there exists a neighborhood $\mcal{S}_0\subset\mcal{S}_a^-$ such that for all $\bm\zeta\in\mcal{S}_0$, $\nabla^2\vP^d(\bm\zeta)\preceq0$. Since the map $\bm x= \GG_a^\inv\ff$ is continuous over $\calS_a$, the image of the map over $\mcal{S}_0$ is a neighborhood of $\bxx$, which we denote as $\mcal{X}_0$. Next, we are going to prove that for any $\bm x\in\mcal{X}_0$, $\nabla^2\vP(\bm x)\preceq0$, which plus the fact that $\bar{\bm x}$ is a critical point of $\vP(\bm x)$ implies $\bar{\bm x}$ is a maximizer of $\vP(\bm x)$ over $\mcal{X}_0$. For any $\bm x\in\mcal{X}_0$, let $\bm\zeta$ be a point satisfying $\bm x= \GG_a^\inv\ff$. Thus, $\nabla^2\vP^d(\bm\zeta)=- \FF^T \GG_a^\inv \FF- \DD^{-1}\preceq0$. By singular value decomposition, there exist orthogonal matrices $ \EE\in\mbb{R}^{n\times n}$, $ \KK\in\mbb{R}^{m\times m}$ and $ \RR\in\mbb{R}^{n\times m}$ with
    \begin{equation}\label{eq:th2R}
     \RR_{ij}=
    \left\{\begin{array}{ll}
    \delta_i, &~~ i=j \wand i=1,\ldots,r,\\
    0,&~~\wotherwise,
    \end{array}\right.
    \end{equation}
    where $\delta_i>0$ for $i=1,\ldots,r$ and $r=\rank( \FF)$, such that
    \begin{equation}\label{eq:th2FDERK}
     \FF \DD^{\frac{1}{2}}= \EE \RR \KK.
    \end{equation}
    Then we have
    \begin{equation}\label{eq:th2hessPid}
     -\DD^{-1}- \DD^{-\frac{1}{2}} \KK^T \RR^T \EE^T \GG_a^\inv \EE \RR \KK \DD^{-\frac{1}{2}}\preceq0.
    \end{equation}
    Being multiplied by $ \KK \DD^{\frac{1}{2}}$ from the left and $ \DD^{\frac{1}{2}} \KK^T$ from the right, this equation can be converted equivalently into
    \begin{equation}\label{eq:th2hessPidequ}
     -\II_m- \RR^T \EE^T \GG_a^\inv \EE \RR\preceq0,
    \end{equation}
    which, by Lemma \ref{lm:PplusDUD} in Appendix, is further equivalent to
    \begin{equation}\label{eq:th2hessPidequfur}
     \EE^T \GG_a \EE+ \RR \RR^T\preceq0.
    \end{equation}
Multiplying the equation (\ref{eq:th2hessPidequfur}) by $ \EE$ from the left and $ \EE^T$ from the right, we obtain
    \begin{equation}\label{eq:th2hessPi}
    0\succeq \GG_a+ \EE \RR \KK \DD^{-\frac{1}{2}} \DD \DD^{-\frac{1}{2}} \KK^T \RR^T \EE^T= \GG_a+ \FF \DD \FF^T=\nabla^2\vP(\bm x).
    \end{equation}
    Therefore, $\bar{\bm x}$ is a maximizer of $\vP(\bm x)$ over $\calX_0$.

    Similarly, we can prove that if $\bar{\bm x}$ is a maximizer of $\vP(\bm x)$ over $\calX_0$, $\bar{\bm\zeta}$ is a maximizer of $\vP^d(\bm\zeta)$ over $\mcal{S}_0$. Then, by the Theorem \ref{th:AnalSolu}, the equation (\ref{eq:th2maxmax}) is proved.

\item We then prove the double-min duality.
\benum
\item  Suppose that $\bar{\bm\zeta}$ is a local minimizer of $\vP^d(\bm\zeta)$ over $\mcal{S}_a^-$. Then there exists a neighborhood $\mcal{S}_0\subset\mcal{S}_a^-$ of $\bar{\bm\zeta}$ such that for any $\bm\zeta\in\mcal{S}_0$, $\nabla^2\vP^d(\bm\zeta)\succeq0$. Let $\mcal{X}_0$ denote the image of the map $\bm x= \GG_a^\inv\ff$ over $\mcal{S}_0$, which is a neighborhood of $\bxx$.
    For any $\bm x\in\mcal{X}_0$, let $\bm\zeta$ be a point that satisfies $\bm x= \GG_a^\inv\ff$. From $\nabla^2\vP^d(\bm\zeta)=- \FF^T \GG_a^\inv \FF- \DD^{-1}\succeq0$, we have $- \FF^T \GG_a^\inv \FF\succeq \DD^{-1}\succ0$, which implies that the matrix $ \FF$ is invertible. Then we obtain
    \begin{equation}\label{eq:th2iii1hessPid}
    - \GG_a^\inv\succeq( \FF^T)^{-1} \DD^{-1} \FF^{-1},
    \end{equation}
  which is further equivalent to
    \begin{equation}\label{eq:th2iii1hessPi}
    - \GG_a\preceq \FF \DD \FF^T.
    \end{equation}
   Hence, we prove that $\nabla^2\vP(\bm x)= \GG_a+ \FF \DD \FF^T\succeq0$ and $\bxx$ is a local minimizer of $\vP(\bm x)$. The converse can be proved similarly. By the Theorem \ref{th:AnalSolu}, the equation (\ref{eq:th2minmin}) is then proved.

\item Suppose that $\bar{\bm\zeta}$ is a local minimizer of $\vP^d(\bm\zeta)$ over $\mcal{S}_a^-$. We claim that $\bar{\bm x}$ is not a local minimizer of $\vP(\bm x)$. If $\bar{\bm x}$ is a local minimizer of $\vP(\bm x)$, we would have $\nabla^2\vP(\bar{\bm x})= \GG_a+ \FF \DD \FF^T\succeq0$, which is equivalent to $ \FF \DD \FF^T\succeq- \GG_a$. Since $- \GG_a\succ0$, it is true that matrix $ \FF$ has full rank and
    \begin{equation}\label{eq:th2iii2rank}
    n=\rank(- \GG_a)=\rank( \FF \DD \FF^T)\leq\min\lbc\rank( \FF),\rank( \DD)\rbc= m,
    \end{equation}
    which is a contradiction. Therefore, plus the previous discussion, $\bar{\bm x}$ must be a saddle point of $\vP(\bm x)$.

    Let
    \begin{equation}\label{eq:th2iii2phi}
    \varphi(\bm t)=\vP(\bar{\bm x}+\LL\bm t).
    \end{equation}
    It can be proved that $\bm0\in\mbb{R}^m$ is a local minimizer of the function $\varphi(\bm t)$, because we have
    \begin{eqnarray}
    &&\nabla\varphi(\bm0)= \LL^T\nabla\vP(\bar{\bm x})=0,\label{eq:th2iii2phipar}\\
    &&\nabla^2\varphi(\bm0)= \LL^T\nabla^2\vP(\bar{\bm x}) \LL\succeq0.
    \end{eqnarray}
    Thus the equation (\ref{eq:th2minmln}) is true.

\item   The proof is similar to that of case (b).
\eenum
\eenum
The theorem is proved. \hfill\qed

\section{Two special problems}\label{se:specific}
In this section, we discuss two special problems, a fourth-order polynomial minimization problem and a minimax problem. The min-max duality is reinforced and existence conditions for the critical point in $\calS_a^+$ are derived. If the existence conditions does not hold, which leads to the nonexistence of critical points in $\calS_a^+$, the perturbation method proposed in \cite{Chen13QPS} can then be applied to solve them approximately.
\subsection{A fourth-order polynomial minimization problem}
The fourth-order polynomial minimization problem considered here is
\begin{equation}\label{p:fourthpoly}
(\calP_1)~~\min_{\xx\in\bbR^n} \vP_1(\xx)=\frac{1}{2}\bm x^T \AA\bm x-\bm f^T\bm x+\frac{\a}{2}\lpa \half\xx^T\BB\xx+c \rpa^2
\end{equation}
where the matrix $\BB$ is symmetric and positive definite.  Without loss of generality, here we assume $\BB=\II$, the identity matrix.

Let
$$
\GG_a=\AA+\s \II, \wand \calS_a^+=\lbc  \s ~|~ \GG_a\succ0, \s\geq\a c \rbc.
$$
The canonical dual problem can be defined as
\begin{equation}\label{p:fourthpolydual}
(\calP_1^d)~~\min_{\s\in\calS_a^+} \vP_1^d(\s)=-\half\ff^T\GG_a^{\inv}\ff+c\s-\frac{1}{2\a}\s^2
\end{equation}

The eigendecomposition of $\AA$ is $\AA=\UU\LL \UU^T$, where diagonal entities of the diagonal matrix $\LL$ are eigenvalues and columns of the matrix $\UU$ are the corresponding eigenvectors. We assume that diagonal entities of $\LL$ are in nondecreasing order, i.e.
$$
\l_1=\cdots=\l_k<\l_{k+1}\leq \cdots\leq \l_{n}.
$$
 If let $\hff=\UU^T\ff$, the dual function can then be rewritten as
\[\label{eq:4polydualhat}
\vP_1^d(\s)=-\half \sum_{i=1}^n\frac{\hatf_i^2}{\l_i+\s}+c\s-\frac{1}{2\a}\s^2
\]

The first-order and second-order derivatives of the dual function $\vP_1^d(\s)$ are
\[
&\d\vP_1^d(\s)=\half \sum_{i=1}^n\frac{\hatf_i^2}{(\l_i+\s)^2}+c-\frac{1}{\a}\s\non\\
&\d^2\vP_1^d(\s)=-\sum_{i=1}^n\frac{\hatf_i^2}{(\l_i+\s)^3}-\frac{1}{\a}\non
\]
Since $\a$ is assumed to be positive, $\d^2\vP_1^d(\s)$ is negative over $\calS_a^+$, which indicates that the dual function is concave over $\calS_a^+$. If $\a c>-\l_1$, $\calS_a^+=[\a c, +\infty)$ and the maximizer of  $\vP_1^d(\s)$ in $\calS_a^+$ corresponds to the unique global solution of the primal problem, since $\GG_a$ is positive definite. If $\a c\leq-\l_1$, $\calS_a^+=(-\l_1, +\infty)$ and we have the following theorem about the existence of a critical point in $\calS_a^+$. Its proof is similar to that in \cite{Chen13QPS}.

\begin{proposition}\label{th:exist_poly}
{\rm (\textbf{Existence Conditions})}
Suppose that $\l_i$ are defined as above and $-\l_1\geq\a c$. Then there exists a critical point of $\vP_1^d(\s)$ in $\calS_a^+$ if and only if  $\sum_{i=1}^k\hatf_i^2\neq0$ or $\half\sum_{i=k+1}^n\hatf_i^2/(\l_i-\l_1)^2+\l_1/\a+c>0$. If $\vP_1^d(\s)$ has a critical point in $\calS_a^+$, the critical point is unique. Let $\bgs$ denote the critical point. Then $\bxx=\GG_a(\bgs)^\inv\ff$ is a global solution of the problem $(\calP_1)$.
\end{proposition}

\subsection{A minimax problem}
In this section, we consider the minimax problem:
\[\label{p:minimax}
\min_{\xx\in\bbR^n}\max\lbc\half\xx^T\AA_1\xx-\ff_1^T\xx+d_1,\half\xx^T\AA_2\xx-\ff_2^T\xx+d_2\rbc
\]
where the matrices $\AA_1$ and $\AA_2$ are symmetric.  Here, we assume that  $\AA_2-\AA_1$ is positive definite.
Generally, the max function is unsmooth, and it can be approximated by the smooth log-sum-exp function
$$
\frac{1}{\b}\log\lpa \exp\lpa \b(\half\xx^T\AA_1\xx-\ff_1^T\xx+d_1) \rpa + \exp\lpa \b(\half\xx^T\AA_2\xx-\ff_2^T\xx+d_2) \rpa\rpa.
$$
The parameter $\b$ is a positive number, and a better approximation will be get if a larger $\b$ is equipped. If we move the function $\half\xx^T\AA_1\xx-\ff_1^T\xx+d_1$ out of the log-sum-exp function, the first item in the log function will become 1 and the quadratic function in the second item will be $\half\xx^T(\AA_2-\AA_1)\xx-(\ff_2-\ff_1)^T\xx+d_2-d_1$. Since $\AA_2-\AA_1$ is assumed to be positive definite, the quadratic function can be transformed into $\half\xx^T\xx+d$ if the coordinate system is properly rotated and moved. Without loss of generality, we consider the following problem:
\[\label{p:Ps}
(\calP_2)~~\min_{\xx\in\bbR^n} \vP_2(\xx)=\frac{1}{2}\bm x^T\AA\bm x-\ff^T\xx+\frac{1}{\beta}\log\lpa 1+\exp\lpa\beta\lpa\frac{1}{2}\bm x^T \xx+d\rpa\rpa\rpa.
\]

The dual function will be an univariate function with
$$
\GG_a=\AA+\t \II.
$$
Here, the $\calS_a^+$ is defined as
\[\label{eq:PsSa}
\calS_a^+=\{\t~|~ 0<\t<1, \GG_a\succ 0\},
\]
and the canonical dual problem is
\[\label{p:Psdual}
(\calP_2^d)~~ \max_{\t\in\calS_a^+}~~\vP_2^d(\t)=-\half \ff^T\GG_a^\inv\ff+d\t-\frac{1}{\b}\lpa \t\log(\t)+(1-\t)\log(1-\t) \rpa
\]

Similarly, the $\AA$ can be decomposed as $\AA=\UU\LL \UU^T$. The diagonal entities of $\LL$ are the eigenvalues of the matrix $\AA$ in nondecreasing order,
$$
\l_1=\cdots=\l_k<\l_{k+1}\leq \cdots\leq \l_{n}.
$$
The columns of $\UU$ are the corresponding eigenvectors. If we let $\hff=\UU^T\ff$, the dual function can be rewritten as
\[\label{eq:Psdualhat}
\vP_2^d(\t)=-\half\sum_{i=1}^n\frac{\hatf_i^2}{\l_i+\t}+d\t-\frac{1}{\b}\lpa \t\log(\t)+(1-\t)\log(1-\t) \rpa.
\]
The set $\calS_a^+$ becomes
$$
\calS_a^+=\{\t~|~ 0<\t<1, \t>-\l_1\}.
$$

It can be noticed that if $\l_1\geq0$, i.e. the matrix $\AA$ is positive semidefinite, the matrix $\GG_a$ is always positive definite as $\t\in\calS_a^+=\{\t~|~ 0<\t<1\}$. There exist a unique critical point of the dual function in $\calS_a^+$.
If $\l_1\leq -1$, the set $\calS_a^+$ will be empty. It can be proved that the minimization problem is not lower bounded. For the case where $-1<\l_1<0$, existence conditions for the dual function having a critical point in $\calS_a^+$ is derived.
\begin{proposition}\label{th:exist_minimax}
{\rm (\textbf{Existence Conditions})}
Suppose that $\l_i$ are defined as above and $-1<\l_1<0$. Then there exists a critical point of $\vP_2^d(\t)$ in $\calS_a^+$ if and only if  $\sum_{i=1}^k\hatf_i^2\neq0$ or $\half\sum_{i=k+1}^n\hatf_i^2/(\l_i-\l_1)^2-\frac{1}{\b}\log(-\l_1/(1+\l_1))+d>0$. If $\vP_2^d(\t)$ has a critical point in $\calS_a^+$, the critical point is unique. Let $\bgt$ denote the critical point. Then $\bxx=\GG_a(\bgt)^\inv\ff$ is a global solution of the problem $(\calP_2)$.
\end{proposition}

\section{Examples}\label{se:exs}

In this section, three examples are provided to illustrate the perfect duality of the canonical duality theory.

\subsection*{Example 1}
\begin{figure}
\begin{center}
\scalebox{0.7}[0.7]{\includegraphics{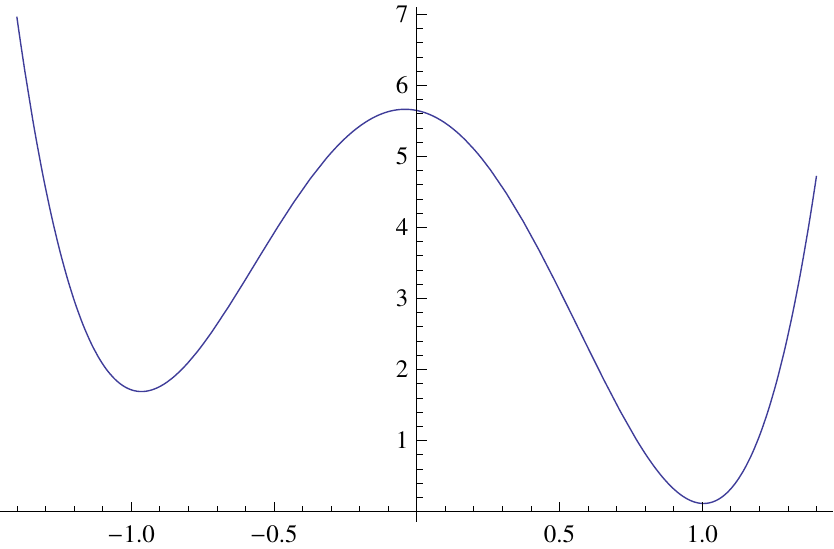}}
\end{center}
\caption{The graph of $\vP(x)$ in Example 1.}
\label{fig:ex3pi}
\end{figure}

\begin{figure}
\begin{center}
\scalebox{0.8}[0.8]{\includegraphics{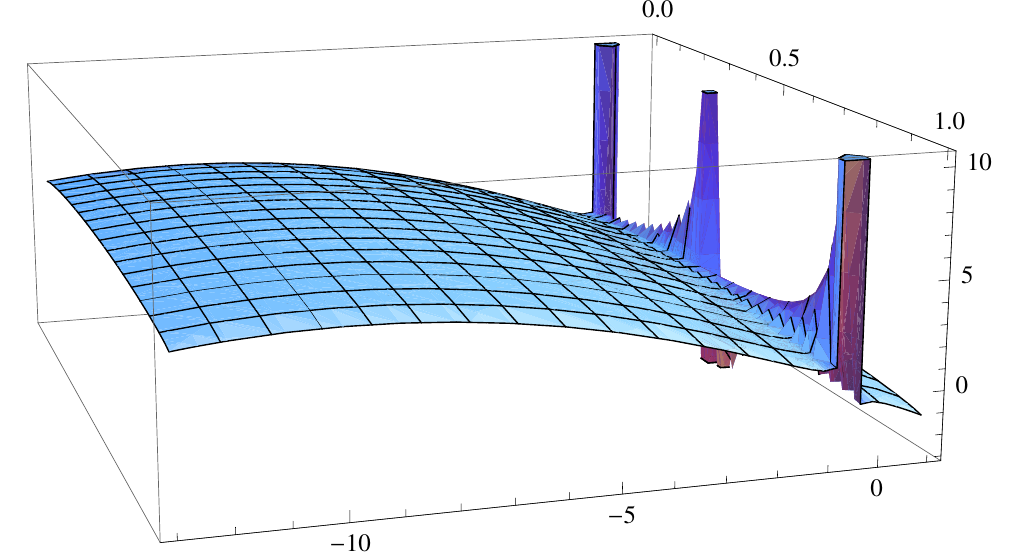}}
\scalebox{0.6}[0.6]{\includegraphics{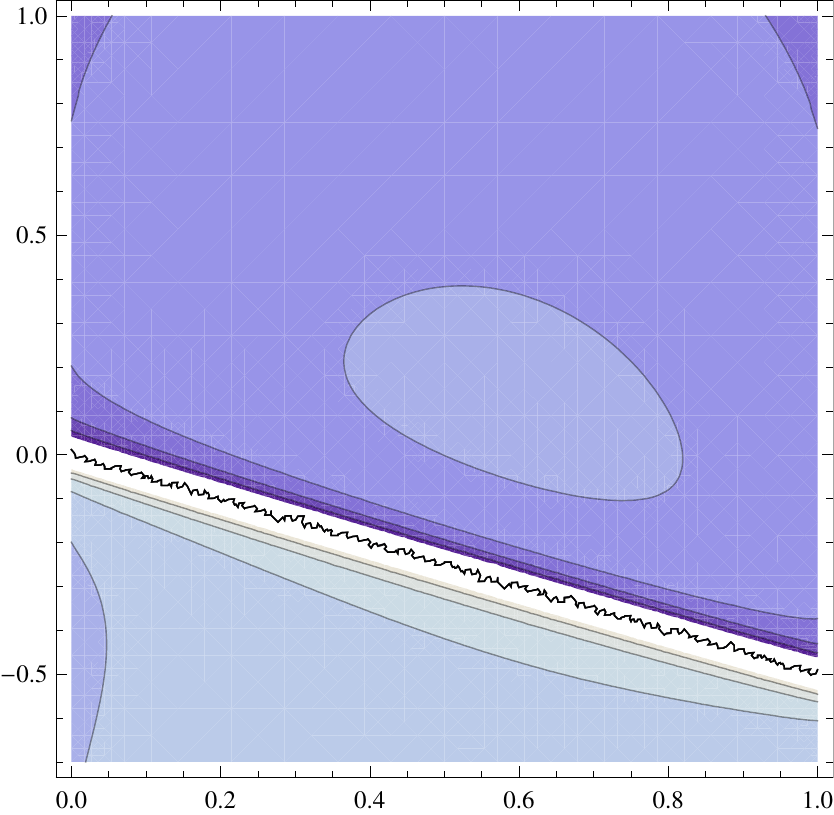}}
\end{center}
\caption{The graph of $\vP^d(\tau,\sigma)$ (left) and  the contour plot of $\vP^d(\tau,\sigma)$ on the region $\calS_a^+$ (right) in Example 1.}
\label{fig:ex3pid}
\end{figure}

~~Consider the 1-dimensional problem:
\[
\min_{x\in\bbR} &\quad\vP(x)=\log\lbk1+\exp\lpa0.5x^2-0.1\rpa\rbk+5\lpa x^2-1\rpa^2-0.8x.\non
\]
The corresponding canonical dual function is
$$
\vP^d(\tau,\sigma)=-\frac{0.32}{\tau+2\sigma}-\sigma-0.05\sigma^2-0.1\tau-\lbk\tau \log(\tau)+(1-\tau)\log(1-\tau)\rbk.
$$
The graph of function $\vP(x)$ is shown in Figure \ref{fig:ex3pi},  and the graph and contour plot of $\vP^d(\tau,\sigma)$ is shown in Figure \ref{fig:ex3pid}.

There are three critical points of the dual function $\vP^d(\tau,\sigma)$:
$$
\begin{pmatrix}\bar{\t}_1\\\bar{\s}_1\end{pmatrix}=\begin{pmatrix}0.599866 \\0.098119 \end{pmatrix},~
\bpmat \bar{\t}_2\\\bar{\s}_2\epmat=\bpmat 0.475231\\-9.983154 \epmat, \wand
\bpmat \bar{\t}_3\\\bar{\s}_3\epmat=\bpmat 0.590128\\-0.71007 \epmat
$$
which are corresponding to the solutions of the primal problem:
$$
\bar{x}_1=1.004894,~\bar{x}_2=-0.041044, \wand \bar{x}_3=-0.963843.
$$
It is noticed that $(\bgt_1,\bgs_1)$ is in $\calS_a^+$ and $\barx_1$ is the global solution of the primal problem, which demonstrates the min-max duality. The double-max duality can be seen from the fact that $(\bgt_2,\bgs_2)$ and $\barx_2$ are local maximizers of functions $\vP(x)$ and $\vP^d(\t,\s)$. Since $n=1$, $ m=2$ and $\barx_3$ is a local minimizer of the function $\vP(x)$, the fact that $(\bgt_3,\bgs_3)\in\calS_a^-$ is a saddle point of the function $\vP^d(\t,\s)$ illustrates the double-min duality.
Moreover, we have
\[
&\vP(\barx_1)=\vP^d(\bgt_1,\bgs_1)=0.112521,\non\\
&\vP(\barx_2)=\vP^d(\bgt_2,\bgs_2)=5.660800,\non\\
&\vP(\barx_3)=\vP^d(\bgt_3,\bgs_3)=1.688196.\non
\]

\subsection*{Example 2}
\begin{figure}
\begin{center}
\scalebox{0.6}[0.6]{\includegraphics{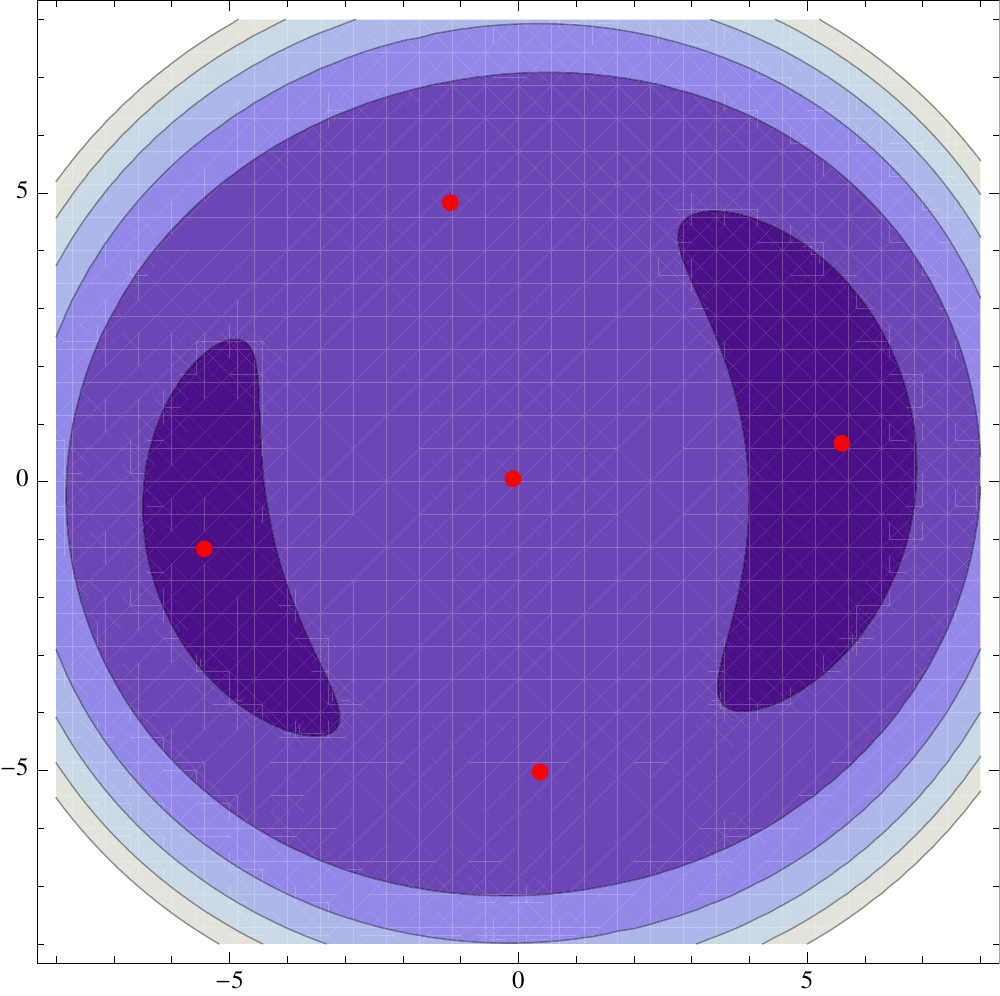}}
\scalebox{0.7}[0.7]{\includegraphics{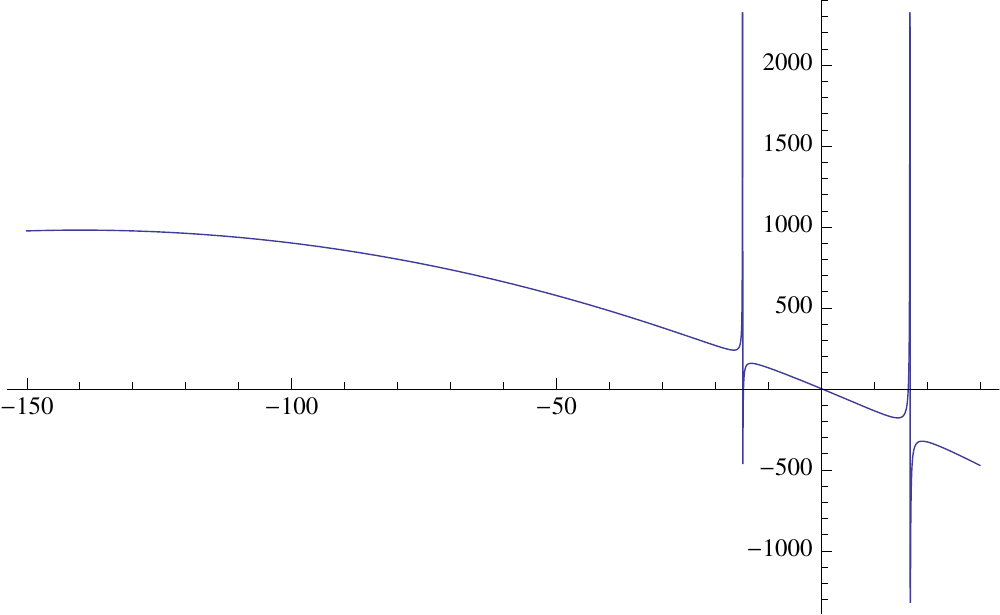}}
\end{center}
\caption{The contour plot of $\vP(\xx)$ and the graph of the dual function in the Example 2.}
\label{fig:ex_poly}
\end{figure}

~~Consider a randomly generated fourth-order polynomial problem:
$$
\min_{\xx\in\bbR^n} \vP(\xx)=\frac{1}{2}\bm x^T \AA\bm x-\bm f^T\bm x+\frac{\a}{2}\lpa \half\xx^T\xx+c \rpa^2
$$
with
$$
\AA=\bpmat -16 & -5 \\ -5 & -14 \epmat,~ \ff=\bpmat 14\\ -6 \epmat, ~ c=-14, \wand \a=10;
$$
The contour plot of $\vP(\xx)$ and the graph of the dual function are shown in Figure \ref{fig:ex_poly}.

There are five critical points of the dual function:
$$
\bgs_1=19.093,  \bgs_2=14.495, \bgs_3=-13.184, \bgs_4=-16.459, \wand \bgs_5=-139.945.
$$
The eigenvalues of the corresponding matrix $\GG_a$ are
$$
\gll_1=\bpmat 2.282\\33.904\epmat, \gll_2=\bpmat -2.32\\29.31\epmat, \gll_3=\bpmat -29.99\\1.63\epmat, \gll_4=\bpmat -33.27\\-1.65\epmat, \wand \gll_5=\bpmat -156.76\\-125.13\epmat,
$$
and the corresponding critical points of the primal problem are:
$$
\xx_1=\bpmat 5.6\\0.67\epmat, \xx_2=\bpmat -5.44\\-1.16\epmat, \xx_3=\bpmat 0.38\\-5.02\epmat, \xx_4=\bpmat -1.18\\4.83\epmat, \wand \xx_5=\bpmat -0.09\\0.05\epmat.
$$

We notice that $\bgs_1$ is in $\calS_a^+$ and $\bxx_1$ is the global solution of the primal problem, which illustrates the min-max duality. Both $\bgs_4$ and $\bgs_5$ are in $\calS_a^-$, the double-min duality is demonstrated by the fact that $\bgs_4$ is a local minimizer  and $\xx_4$ is a saddle point, and the double-max duality is demonstrated by the fact that $\bgs_5$ is a local maximizer and $\xx_5$ is a local maximizer. Moreover, the values of the primal function and dual function are equal on each pair of solutions.

\subsection*{Example 3 {\rm (\cite{kiwiel85book})}}
\begin{figure}
\begin{center}
\scalebox{0.7}[0.7]{\includegraphics{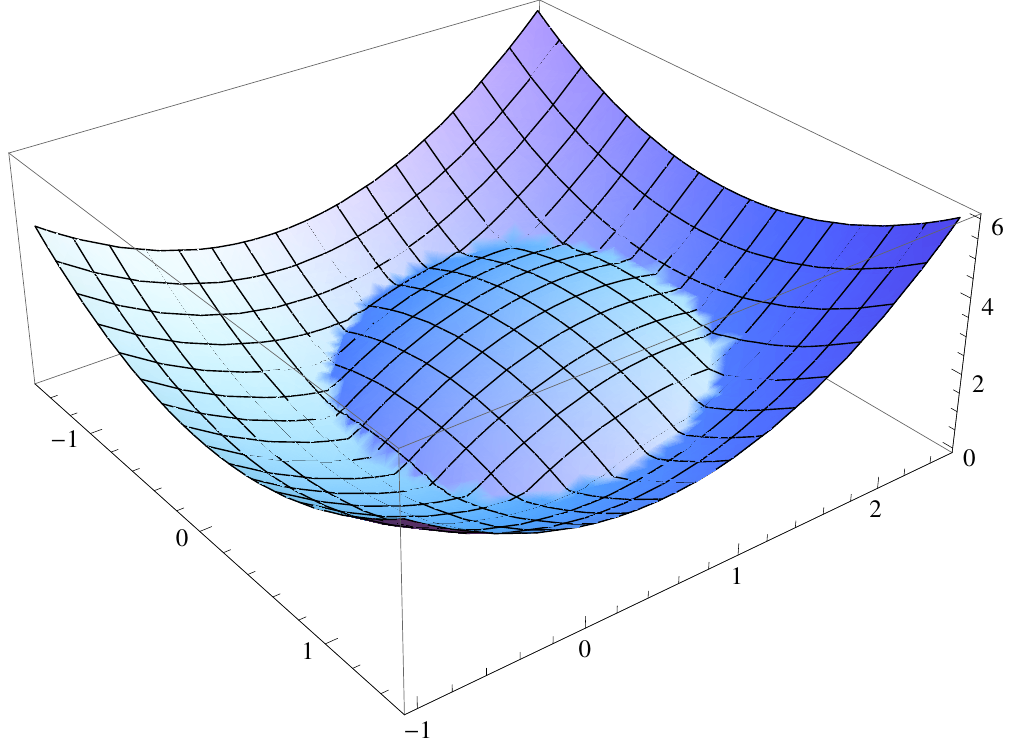}}
\scalebox{0.7}[0.7]{\includegraphics{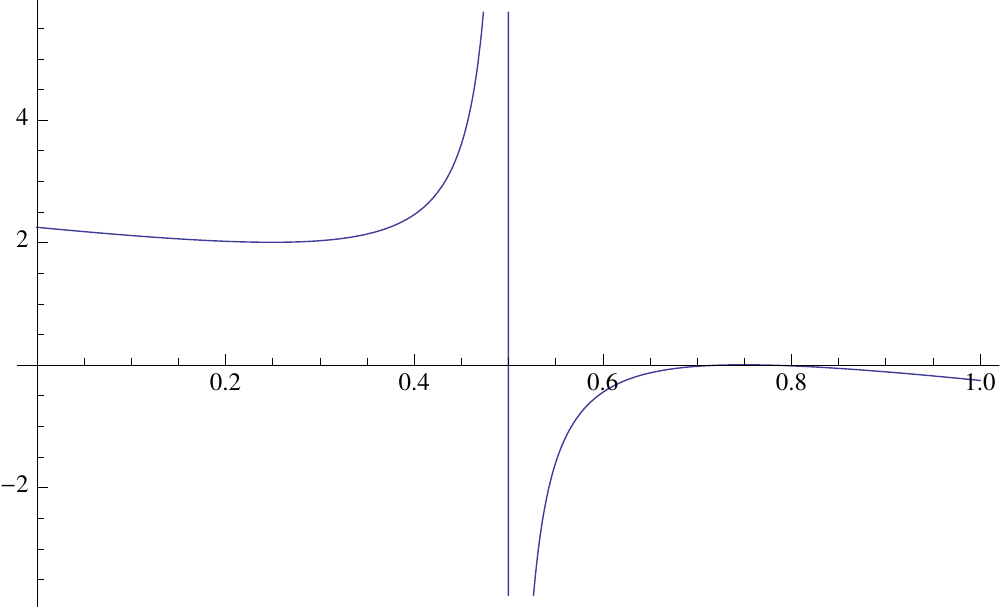}}
\end{center}
\caption{The graph of $\vP(\xx)$ (left) and the graph of $\vP^d(\t)$ (right) in Example 3. The parameter is set as $\b=100$. }
\label{fig:exminimax}
\end{figure}

We consider a nonconvex and nonsmooth optimization problem:
$$
\min_{\xx\in\bbR^2} \max\lbc x_1^2+x_2^2-x_2, ~-x_1^2-x_2^2+3x_2 \rbc.
$$
It's easy to verify that the optimal solution is $(0,0)$ with value 0. Here, we use the log-sum-exp function to approximate the function $\max\lbc\cdot,\cdot\rbc$, and then get the following smooth optimization problem:
$$
\min_{\xx\in\bbR^2} \vP(\xx)=\frac{1}{\b}\log\lbk 1+ \exp\lpa \b\lpa 2x_1^2+2x_2^2-4x_2\rpa \rpa\rbk-x_1^2-x_2^2+3x_2.
$$
Its canonical dual function is
$$
\vP^d(\t)=-\half\frac{(4\t-3)^2}{4\t-2}-\frac{1}{\b}\lbk\t\log\t+(1-\t)\log(1-\t)\rbk.
$$
The graphs of the approximation function $\vP(\xx)$ and the dual function $\vP^d(\t)$ are shown in Figure \ref{fig:exminimax}.

With $\b=100$, the two critical points of the dual function $\vP^d(\t)$ are
$$
\bgt_1=0.749318, \wand \bgt_2=0.249308.
$$
The corresponding solutions of the primal problem are
$$
\bxx_1=\bpmat 0\\-0.002734\epmat, \wand \bxx_2=\bpmat0\\1.99724\epmat.
$$
The min-max duality is true by the fact that $\t_1\in\calS_a^+$ and $\xx_1$ is the global solution of the primal problem. The double-min duality is also true because of the fact that $ m=1$, $n=2$, $\t_2\in\calS_a^-$ is a local minimizer and $\xx_2$ is a saddle point of the primal function. Moreover, we have
$$
\vP(\bxx_1)=\vP^d(\bgt_1)=0.005627, \wand\vP(\bxx_2)=\vP^d(\bgt_2)=2.00562.
$$

\section{Conclutions}\label{se:concl}
A very general nonlinear global optimization problem with a 4th-order polynomial and a log-sum-exp function is discussed. The canonical duality theory is applied to deal with this challenging problem. The triality theory concludes that if there is a critical point in the positive semidefinite region $\calS_a^+$, the dual problem can be solved easily and, correspondingly, the global solution of the primal problem can be found analytically from this critical point. For the two specific problems, a fourth-order polynomial minimization problem and a minimax problem, existence conditions for the critical point are derived. If these conditions hold, there must be a critical point in the positive semidefinite region and then the global solution of the primal problem can be easily obtained by solving the dual problem. The examples demonstrate the perfect duality of the canonical duality theory.

Besides the duality about the global solution, which is called the min-max duality, the triality theory also discusses the relationships of local extrema, which are called double-min duality and double-max duality. But the saddle points which are not in the $\calS_a^+$ and $\calS_a^-$ have not been clarified. There are some interesting phenomenas, which hint that these saddle points may can be sorted according to a certain order.  It is our next work to investigate the order of the extrema in the dual space. Another work would be on constructing existence conditions of critical points for the general problem.

\bibliography{logexp-abbrv}
\bibliographystyle{jogo}

\section*{Appendix}\label{se:appendix}

The following lemma is a generalization of Lemma 6 in \cite{Gao-Wu-triality12}.
\begin{lemma}\label{lm:PplusDUD}
Suppose that $ \PP\in\mbb{R}^{n\times n}$, $ \UU\in\mbb{R}^{m\times m}$ and $ \DD\in\mbb{R}^{n\times m}$ are given symmetric matrices with
$$
 \PP=\begin{bmatrix} \PP_{11} &  \PP_{12}\\ \PP_{21} &  \PP_{22} \end{bmatrix}\prec 0, ~~
 \UU=\begin{bmatrix} \UU_{11} & \bm 0\\\bm 0 &  \UU_{22} \end{bmatrix}\succ 0, \wand
 \DD=\begin{bmatrix} \DD_{11} & \bm 0\\\bm 0 & \bm 0 \end{bmatrix},
$$
where $ \PP_{11}$, $ \UU_{11}$ and $ \DD_{11}$ are $r\times r$-dimensional matrices, and $ \DD_{11}$ is nonsingular. Then,
\begin{equation}\label{eq:apl m}
 \PP+ \DD \UU \DD^T\preceq0 \Leftrightarrow - \DD^T \PP^{-1} \DD- \UU^{-1}\preceq0.
\end{equation}
\end{lemma}

\proof Obviously, $ \PP+ \DD \UU \DD^T\preceq0$ is equivalent to
\begin{equation}\label{eq:aplm-2}
- \PP- \DD \UU \DD^T=
\begin{bmatrix}
- \PP_{11}- \DD_{11} \UU_{11} \DD_{11}^T & - \PP_{12}\\
- \PP_{21} & - \PP_{22}
\end{bmatrix}\succeq0.
\end{equation}
By Schur lemma, equation (\ref{eq:aplm-2}) is equivalent to
\begin{equation}\label{eq:aplm-3}
- \PP_{11}- \DD_{11} \UU_{11} \DD_{11}^T+ \PP_{12} \PP_{22}^{-1} \PP_{21}\succeq0 \wand \PP_{22}\prec0.
\end{equation}
The inverse of matrix $ \PP$ is
$$
 \PP^{-1}=
\begin{bmatrix}
( \PP_{11}- \PP_{12} \PP_{22}^{-1} \PP_{21})^{-1} & - \PP_{11}^{-1} \PP_{12}( \PP_{22}- \PP_{21} \PP_{11}^{-1} \PP_{12})^{-1}\\
-( \PP_{22}- \PP_{21} \PP_{11}^{-1} \PP_{12})^{-1} \PP_{21} \PP_{11}^{-1} & ( \PP_{22}- \PP_{21} \PP_{11}^{-1} \PP_{12})^{-1}
\end{bmatrix}.
$$
Then, it is easy to prove that $- \PP_{11}+ \PP_{12} \PP_{22}^{-1} \PP_{21}\succ0$. Since $ \DD_{11}$ is nonsingular and $ \UU_{11}\succ0$, we have $ \DD_{11} \UU_{11} \DD_{11}^T\succ0$. Thus the equation (\ref{eq:aplm-3}) is equivalent to
\begin{equation}\label{eq:aplm-4}
(- \PP_{11}+ \PP_{12} \PP_{22}^{-1} \PP_{21})^{-1}\preceq( \DD_{11} \UU_{11} \DD_{11}^T)^{-1},
\end{equation}
which is further equivalent to
\begin{equation}\label{eq:aplm-5}
 \DD_{11}^T(- \PP_{11}+ \PP_{12} \PP_{22}^{-1} \PP_{21})^{-1} \DD_{11}\preceq \UU_{11}^{-1}.
\end{equation}
Since $ \DD_{11}^T(- \PP_{11}+ \PP_{12} \PP_{22}^{-1} \PP_{21})^{-1} \DD_{11}=- \DD^T \PP^{-1} \DD$ and $ \UU_{22}\succ0$, the equation (\ref{eq:aplm-5}) is equivalent to
\begin{equation}\label{eq:aplm-6}
- \DD^T \PP^{-1} \DD- \UU^{-1}\preceq0.
\end{equation}
The lemma is proved. \hfill \qed

\end{document}